\documentclass[11pt]{article}
\usepackage{latexsym,amssymb,amsmath}

\newcommand{\ep}{\varepsilon}
\newcommand{\toro}{{\mathbb T ^2}}
\newcommand{\rr}{{\mathbb R}}
\newcommand{\zz}{{\mathbb Z ^2}}
\newcommand{\zo}{{\mathbb Z} ^2_0}
\newcommand{\pa}{\partial}
\newcommand{\na}{\nabla}
\newcommand{\sscri}{\scriptscriptstyle}

\newcommand{\beq}{\begin{equation}}
\newcommand{\eeq}{\end{equation}}
\newcommand{\bsquare}{\hbox{\rule{6pt}{6pt}}}

\newtheorem{proposition}{Proposition}[section]
\newtheorem{theorem}[proposition]{Theorem}
\newtheorem{lemma}[proposition]{Lemma}

\begin{document}

\title{2D Navier-Stokes equation in Besov spaces of negative order}
\author{Zdzis{\l}aw Brze\'zniak \\Department of Mathematics\\
                         The University of York\\
                         Heslington, York YO10 5DD, UK
       \and  Benedetta Ferrario \\Dip. di Matematica ``F. Casorati''\\
                           Universit\`a di Pavia\\
                           I-27100 Pavia}
\date{ }
\maketitle

\begin{abstract}
      The Navier--Stokes equation in the bidimensional torus is considered,
     with initial velocity in the Besov spaces
     $B^{-s+2-\frac{2}{r}}_{p\,r}$ and
     forcing term in $L^r(0,T;B^{-s}_{p\,q})$
     for suitable indices $s,r,p,q$.
     Results of local existence and uniqueness are proven in the case
     $-1<-s+2-\tfrac 2r<0$ and of global existence in the case
     $-\frac 12 <-s+2-\tfrac 2r<0$.
 \end{abstract}
\noindent
{\bf Key words}:
Navier-Stokes equations,
weak solutions,
existence uniqueness and regularity theory.
\\
{\bf MSC2000}:
35Q30, 
76D05, 
76D03, 
35D05, 
35D10. 

\section{Introduction} \label{intro}
Analysis of classical or weak solutions to the Navier-Stokes
equation in a two-dimensional bounded domain has been widely
considered in the literature. A lot of work has been done
on solutions with finite energy, because of their physical
meaning. However, infinite energy solutions are important for
other reasons. On the  one hand, analysis of self-similar
solutions in the whole space involves velocity fields with
infinite energy because of the no rapid enough decreasing to
infinity (see, e.g., \cite{ca,gp} and references therein). On the
other hand, invariant (equilibrium) Gibbs measures are known for
two-dimensional hydrodynamics in bounded domains (see, e.g.,
\cite{afhk,ahk} for the deterministic Euler equation and
\cite{dpd,af} for a stochastic Navier--Stokes equation, whereas
\cite{ac} deals with both the viscous and inviscous problems).
The velocity fields with finite energy are negligible with respect
to these invariant measures. This makes interesting the analysis
of  infinite energy solutions in a bidimensional bounded  spatial
domain.

The aim of this paper is to investigate the Navier--Stokes
evolution problem in a bounded domain $\mathbb D$ of $\rr^2$
with initial velocity  and forcing term of low regularity.
The force will be integrable in time to some power with value in some
Besov space
$B^{\sigma}_{p\,q}$ of negative
order $\sigma$ (for the space variable) and the initial velocity
will belong to some other Besov space of negative order.
We point out the importance of
the force; for instance, this allows to consider problems related to the
Navier--Stokes equation with a stochastic forcing term.
(See, e.g., \cite{dpd} for the study of auxiliary (deterministic)
equations of Navier--Stokes type  arising from a stochastic
Navier--Stokes equation. There the forcing
term has space regularity  of a negative order Besov space.
We also refer to the bibliography of \cite{dpd}
for other papers on the stochastic
Navier--Stokes equation.)
Therefore  the regularity of the initial velocity is important as well as
the time-space regularity of the driving force.
We shall provide existence and
uniqueness results (global in time), as a generalization of the
classical results in Hilbert spaces
(recalled in the remark at the end of Section \ref{f-smooth}).

The content of the paper is as follows. In the next section, we shall
introduce the Navier--Stokes equation and the Besov spaces to work with.
In Section \ref{sec:locale}, we shall deal with  existence and
uniqueness results
of solutions to our problem on a small time interval, when the initial
velocity belongs to a Besov space of negative order and the forcing
term is integrable in time to some power
with value into a Besov space of negative order
in space (see Theorem \ref{t1}).
In Section \ref{sec:globale} global (in time) existence results
will be given; for this aim,
first we shall split our problem into two auxiliary problems
(considering the additive splitting of the velocity field
$u=x+y$). The equation for the variable $y$ will have
small initial data and small forcing term; in Proposition \ref{y-globale}
existence of a unique solution $y$ will be proven.
The problem for the variable $x$ will be solved in Proposition
\ref{x-globale} using a priori energy estimates,
 as in  the classical case discussed e.g. in \cite{temam}.
Theorem \ref{teo-gen} will concern the main result
for the given problem in the variable $u$.
Finally, we shall make some remarks on the case of a smooth forcing term
in Section \ref{f-smooth}.
The range of variability of the parameters involved in defining the functional
spaces will be specified in Appendix A and B, providing examples
of non classical solutions.
\section{The Navier-Stokes equation}
We consider the Navier--Stokes evolution problem,
i.e. the system of equations governing the motion of an homogeneous
 incompressible viscous fluid
\begin{equation}\label{pde}
\left\{
  \begin{array}{rl} \displaystyle
     \frac{\pa u}{\pa t}-\nu \Delta u +(u\cdot\nabla)u+\nabla p
                          & =\varphi   \\[2mm]
     \nabla \cdot u&=0
  \end{array}
\right.
\end{equation}
where the spatial domain is the torus
$\toro=\mathbb{R}^2/(2\pi \mathbb{Z})^2$. In other words, we
consider our problem on the square $[0,2\pi]^2$ with  periodic
boundary conditions. We consider a finite time interval $[0,T]$.
Here the unknowns are the velocity vector field
$u=(u_1(t,\xi),u_2(t,\xi))$ and the pressure field of the fluid
$p=p(t,\xi)$ (for $(t,\xi)\in [0,T]\times \toro$); $ \Delta$ is
the Laplacian operator $\Delta=\pa^2/\pa \xi^2_1+\pa^2/\pa
\xi_2^2$, $\na$ is the gradient operator
$\na=(\pa/\pa \xi_1, \pa/ \pa \xi_2)$ and $\cdot$ is the scalar
product in $\rr^2$. $\varphi$ is a given external force and
$\nu>0$ is the (constant) viscosity coefficient. $\na\cdot u=0$ is
the incompressibility condition.

The basic functional spaces  to set the problem are the
generalized Sobolev spaces
     $H^s_p$ ($s\in \mathbb R, 1 < p <\infty$)
defined to be  ``regular'' subspaces of the space $\mathcal
D^\prime$ of periodic divergence-free vector distributions. To
define them, we proceed in the following way. Take any
$u\in \mathcal{D}^\prime$. Then since  $\nabla\cdot u=0$
on $\toro$, there exists a periodic scalar distribution
$\psi$ on $\toro$, called the stream function, such
that
\begin{equation}\label{stream}
 u=\nabla^\perp \psi
   \equiv (-\pa \psi/\pa \xi_2, \pa \psi/\pa \xi_1).
\end{equation}
Decomposing $\psi$ in Fourier series with respect to the complete
orthonormal system in $L_2(\toro)$
given by $\{ \frac 1{2\pi}  e^{i k \cdot \xi} \}_{k \in \zz} $
$$
 \psi(\xi)=\sum_{k \in \zz} \psi_k \frac{e^{i k \cdot \xi}}{2\pi}
$$
by \eqref{stream} we get that $u$ has the following
 Fourier series representation
\begin{equation}\label{Fou-vel}
 u(\xi)=\sum_{k \in \zo} u_k e_k(\xi),
         \qquad
 u_k \in \mathbb C, \quad \overline u_k = -u_{-k},
\end{equation}
where $e_k(\xi)=\frac{k^\perp }{2 \pi |k|}\,e^{i k \cdot
\xi}$. Note that $\{e_k\}_{k \in \zo}$ is a  complete orthonormal
system of the eigenfunctions (with corresponding eigenvalues
$|k|^2$) of the operator $-\Delta$ in
$[L^{\textrm{div}}_2(\toro)]^2=\{u \in [L_2(\toro)]^2: 
\,\int_\toro u(\xi)\ d\xi=0, \na \cdot
u=0,$ with  the normal component of $u \text{
being periodic  on } \pa \toro\;\}$, $k^\perp=(-k_2,k_1)$,
$|k|=\sqrt{k_1^2+k_2^2}$ and $\zo=\{k \in \zz:  |k|\neq 0 \}$.
\\
Each $e_k$ is a periodic divergence-free $C^\infty$-vector
function (i.e. $e_k \in \mathcal D$). The convergence of the
series \eqref{Fou-vel} depends on the regularity of the
vector function $u$, and can be used to define Sobolev spaces as
in the following definition.\\
For any $s \in \rr, 1<p<\infty$, we define 
\beq \label{def:h}
 H^s_p=\big\{\; u\in \mathcal{D}^\prime:  u
 \underset{\mathcal D^\prime}{=} \sum_{k \in \zo} u_k e_k
 \text{ and }
           \sum_{k \in \zo} u_k |k|^s e_k \in [L_p(\toro)]^2\; \big\}
\eeq 
$H^s_p$ is a Banach space with norm $\|u\|_{H^s_p}=\|\sum_k
u_k |k|^s e_k \|_{L_p}$.
\\
Thus the unknown velocity will be considered as a function of the time
variable taking values in some $H^s_p$ space
$$
 u(t)=\sum_k u_k(t) e_k
$$
\\
Let  $P$ be the following projection  from the space of
periodic distributions onto the space of periodic
divergence-free distributions: $P u= \sum_k \langle
u,e_{-k}\rangle e_k$. Here  $\langle\cdot,\cdot\rangle$ is
the $([C^\infty(\toro)]^2)^\prime- [C^\infty(\toro)]^2$ duality
bracket. Applying the projection $P$ to the first equation
\eqref{pde} we get rid of the pressure term (because $\langle
\nabla p,e_{-k}\rangle=-\langle p, \nabla \cdot e_{-k}\rangle=0$
for any $k$), obtaining the following formulation of  our problem
\begin{equation} \label{ns}
 u^\prime(t)+Au(t)+B(u(t))=f(t),  \qquad t \in (0,T],
\end{equation}
as an equality in the distributional sense, with some
initial condition $u(0)$ assigned.
We have taken the viscosity $\nu=1$, without loss of generality.
\\
$A$ is the stokes operator $A=-\Delta$,
$B$ is the quadratic operator $B(u)=B(u,u)$ defined by the bilinear
operator $B(u,v)=P[(u\cdot \nabla)v]$. Notice that  $B(u,v)$ is
equal to $P[\nabla \cdot (u\otimes v)]$ because of the divergence-free
condition. 
 $f=P\varphi$. 

We are interested in  the evolution problem \eqref{ns} for initial
data and forcing term not too regular (in space variables). For
this purpose, we define a scale of spaces $B^s_{p\,q}$
consisting of the Besov spaces of periodic divergence-free vector
fields. They can be introduced as real interpolation spaces (see,
e.g. \cite{bl} Theorem 6.4.5): \beq \label{def:b}
 B^{s^*}_{p\,q}=(H^{s_1}_p,H^{s_2}_p)_{\theta,q}
     \qquad
 s^*=(1-\theta)s_1+\theta s_2, \quad 0<\theta<1
\eeq
for $s^* \in \rr, 1 <p,q<\infty$.
\\
In particular, $B_{2\,2}^{s^*}=H^{s^*}_2$ and
$B^0_{2\,2}=[L^{\textrm{div}}_2(\toro)]^2$.

Since $\{e_k\}_{k \in \zo}$ are the eigenvectors of the Stokes
operator, with corresponding eigenvalues $\lambda_k=|k|^2$, then
we can naturally extend $A$ to the whole space
$\mathcal{D}^\prime$ by the following formula: $Au= \sum_k u_k
|k|^2 e_k $ for any $u \in \mathcal D^\prime$ and by means
of this representation it is straightforward to show that the
Stokes operator $A$ is a linear operator in $H^s_p$ with domain
$H^{s+2}_p$; moreover it is a bijective unitary operator
 from $H^{s+2}_p$ onto
$H^s_p$ for any index $s \in \rr$, $1<p<\infty$, and also
from
$B^{s+2}_{p\,q}$ onto $B^s_{p\,q}$ (for any $s \in \rr, 1<p,q < \infty$)
(see, e.g., \cite{l} Theorem 1.1.6 for getting linear operators in
interpolation spaces).
In particular the inverse operator $A^{-1}$ is a linear bounded
operator in each space $B^s_{p\,q}$.
Finally, the operator $A$ generates an analytic semigroup in each $B^s_{p\,q}$.
(We refer, e.g., to \cite{te-p} for classical analysis
of the Navier--Stokes equation on the torus.)
\\
The properties of the Stokes operator are the basis to analyze
equation \eqref{ns} as a perturbation of the linear Stokes problem by
the nonlinear operator $B$. To this end,
the key point is to estimate the operator $B$ in Besov spaces,
as done, e.g., in \cite{ch}.
\\[2mm]
{\tt Remark.}
We believe that all what follows holds if the spatial domain is any bidimensional
smooth bounded domain $\mathbb D$ \footnote{
       For instance, there exists a
       complete orthonormal system $\{e_n\}_{n \in \mathbb N}$
       of eigenvectors of the Stokes operator in the
       space of square integrable divergence-free vector fields
       satisfying homogeneous Dirichlet
       boundary condition, with associated eigenvalues
       $0<\lambda_1\leq \lambda_2\leq \ldots$, $\lambda_n \sim n$ as
       $n \to \infty$
       (see, e.g., \cite{temam}). This allows
       to analyze the linear Stokes operator as presented above for
       the case on the torus.
       }, but we postpone this analysis to a future work.
However, we point out that only in the space periodic case
the expression of the Gibbs  measure constructed by means of the
enstrophy (given in \cite{afhk}) is an invariant measure
for a stochastic Navier--Stokes  evolution
problem. This depends on the fact that
this Gibbs measure is invariant for the Euler flow both on the torus
(see \cite{afhk})
or in any smooth bounded domain
(see \cite{ahk}), but the Euler and  Navier--Stokes boundary
conditions are the same only for $\mathbb D=\toro$.
(For an overview on this subject
see e.g. \cite{af-sw}.)
For this reason and because of the interest described in the
introduction,
 we present our problem considering the spatial domain $\toro$.
We point out that our technique does not apply in the case
the spatial domain is the whole space,
because the inverse of the Stokes operator
is not a bounded operator
in $B^s_{p\,q}(\rr^2)$.
Hence, the techniques  are quite different from
the case of the spatial domain $\mathbb D=\rr^2$.
Moreover, also the assumptions for getting existence and uniqueness
results are different from the case $\mathbb D=\rr^2$; let us consider
\cite{bg}, which is as far as we are aware of
the only paper including a forcing term in the analysis of very rough
solutions of the Navier--Stokes problem.
Biagioni and Gramchev assume
$F \in L^1_{loc}(\rr_+;L^q(\rr^2)), q \in [\frac 43,2[$,
  or $F \in L^1_{loc}(\rr_+;L^1(\rr^2)\cap L^\infty(\rr^2) )$
(setting $F=\nabla^\perp \cdot f$, for $f$ the forcing term in
\eqref{ns}).
\\
Other results for $\mathbb D=\rr^3$ are presented in \cite{CP}.
\hfill $\Box$

\section{Local existence and uniqueness} \label{sec:locale}
In this section we prove a result on uniqueness and  local existence
for equation \eqref{ns} in Besov spaces.
 We use a technique from \cite{b91} based on a theorem of local
diffeomorphism of \cite{vf} (\S 1.1).
\\
Let $u$ be a vector valued function defined on the time interval
$[0,T]$ for $T < \infty$ or on $[0,\infty)$ if
$T=\infty$. Keeping in mind the notation in equation
\eqref{ns}, consider the mapping $\Phi \ni  u \mapsto
\big(u^\prime+Au+B(u),u(0)\big) \in \mathcal E \to \mathcal F, $
for some Banach spaces $\mathcal E$ and $ \mathcal
F=\mathcal F^{\sscri 1}\times \mathcal F^{\sscri 0}$.\footnote{In
practice, these
  spaces will be constructed from $L^\alpha(0,T;B^s_{p\,q})$ and
  $B^s_{p\,q} $ spaces. From now on, $B^s_{p\,q}$ denotes the space for
  vector valued functions defined on the torus, unless otherwise
  specified. By
  the way, we remark that in the notation for the time-integrability, the
  {\it sup}-index has been used, whereas in the notation for the
  space-integrability the {\it inf}-index has been used.}
The Fr\'echet derivative  $d_a \Phi$
 at the point $a=0$
is given by the linear operator $d_0 \Phi: \mathcal E \to \mathcal F,
 u\mapsto \big( u^\prime+Au,u(0) \big)$.
Under the assumption that the quadratic operator $B$
from $\mathcal E$ to $\mathcal F^{\sscri 1}$  is
bounded\footnote{
   The operator $B:\mathcal E \to \mathcal F^{\sscri 1}$
   is bounded if
   $$
     \|B\|:=\sup_{\|u\|_{\mathcal E}\leq 1}\|B(u)\|_{\mathcal F^{\sscri 1}}<\infty.
   $$
   }
and that the linear operator $d_0 \Phi$ is an isomorphism,
Vishik and Fursikov \cite{vf} observe that $\Phi$ is
analytic in a neighbourhood of   $0 \in \mathcal E$ and locally
$\Phi$ has an inverse operator $\Psi$, which is analytic in the
neighbourhood $\mathcal{B}(\sigma_\infty)$ of $\Phi(0) \in
\mathcal F$ of radius $\sigma_\infty = \big( 4 \|(d_0
\Phi)^{{\sscri{- 1}}}\|^2 \|B\| \big)^{-1}$.  This
provides \underline{global} existence of a unique solution $u \in
\mathcal E$ for small forcing term and small initial data in
$\mathcal F$. This solution has analytical dependence on
the initial value and forcing term. We will use this technique in
Proposition \ref{y-globale}.
\\
Moreover, starting from this approach, \cite{b91} shows a result
of \underline{local} existence for any data in $\mathcal F$
by showing that the Vishik-Fursikov procedure works also for
$T<\infty$ and by finding an estimate on $ \sigma_T$.
(The sub-index
$T$ reminds that we work on the time interval
$[0,T]$. The important point is that also the norms in
$\mathcal E$ and in $\mathcal F$ will depend on $T$.)
Roughly
speaking, equation \eqref{ns} is seen as a perturbation by a small
nonlinear term $B$ of the linear equation (which is well posed in
the Hadamard's sense, see Proposition \ref{isom}). The
``smallness'' of $B$ is obtained choosing the time interval small
enough. More precisely, let the time variable vary in the finite
interval $[0,T]$; assume that the quadratic term $B:\mathcal E \to
\mathcal F^{\sscri 1} $ has a norm bounded by
$$
 \|B\| \leq \;C_1 \;T\,^\ep
$$
for some positive constants $C_1$ and $\ep$; then the solution $u$ to the
Navier--Stokes equation \eqref{ns} exists locally in time for any
forcing term in $\mathcal F^{\sscri 1}$
and initial velocity in $\mathcal F^{\sscri 0}$.
Indeed, given $\big(f,u_0\big) \in  \mathcal F $,
let $0<\overline T \leq T$ be such that
\begin{equation} \label{picc}
 \|\big(f,u_0\big)\|_{_{\mathcal F}} \; \overline T\,^\ep \;
< \;
 \frac{1}{4\;  \|(d_0 \Phi)^{\sscri{-1}}\|^2 \;C_1}
\end{equation}
Then  $\big( f,u_0 \big) \in \mathcal B(\sigma_{\overline T})$,
the $\mathcal F$-ball of
radius $\sigma_{\overline T}$; according to
Vishik and Fursikov's result quoted above,  we conclude that there exists  a
unique solution  defined on the time interval $[0,\overline T]$.

These are quite general results. The crucial point is the choice of
the spaces $\mathcal E$ and $\mathcal F$.
In \cite{b91}
regular Hilbert spaces are considered:
$\mathcal E= \{ u \in L^2(0,T;H^2_2):   u^\prime \in L^2(0,T;H^0_2)
\}$ and
$\mathcal F =L^2(0,T;H^0_2)\times H^1_2$.
Here we consider  $\mathcal E= \{ u \in L^r(0,T;B^{-s+2}_{p\,q}): u^\prime \in
L^r(0,T;B^{-s}_{p\,q})\}$ and
$\mathcal F= L^r(0,T;B^{-s}_{p\,q}) \times B^{-s+2-\frac{2}{r}}_{p\,r}$ with
$1<p,q,r<\infty$ and $s \in \rr$.
The norms are defined as
\vspace{3mm}

$
\begin{array}{l}
  \|\big( f,u_0 \big)\|_{_{\mathcal F} }   =
  ( \int_0^T \|f(t)\|^r_{B^{-s}_{p\,q}} dt)^{1/r} + \|u_0\|_{B^{-s+2-\frac{2}{r}}_{p\,r}}
\\[3mm]
  \|u\|_{_{\mathcal E}} =
      (\int_0^T \|u(t)\|^r_{B^{-s+2}_{p\,q}} dt\;+\int_0^T \|u^\prime(t)\|^r_{B^{-s}_{p\,q}} dt)^{1/r}\\[2mm]
\end{array}
$
\\
We are especially interested in the case $-s+2-\frac{2}{r}<0$.
\\[1mm]
In this framework, the main result in this section is as follows.
\begin{theorem} \label{t1}
Let
\\
$\hspace*{1cm}\mathcal E= \{ u \in L^r(0,T;B^{-s+2}_{p\,q}): u^\prime \in
L^r(0,T;B^{-s}_{p\,q})\}$
\\
$\hspace*{1cm}\mathcal F^{\sscri 1}= L^r(0,T;B^{-s}_{p\,q}),\qquad
 \mathcal F^{\sscri 0}=B^{-s+2-\frac{2}{r}}_{p\,r},\qquad
 \mathcal F=\mathcal F^{\sscri 1} \times \mathcal F^{\sscri 0}$
\\
with $1<p,q,r<\infty,  \,s \in \mathbb R$.\\
Then
\begin{itemize}
\item[i)]
 The map $d_0 \Phi:\mathcal E \ni u \mapsto \big( u^\prime+Au,u(0)\big) \in \mathcal F$
is well defined and continuous. Moreover, the operator $d_0 \Phi$
is an isomorphism of Banach spaces $\mathcal E$ and $\mathcal F$.
\item[ii)] If the parameters $r,p,q,s$ satisfy the following
conditions
\begin{eqnarray}
  r  \leq q                           \label{1} \\
  s -1< \tfrac{2}{p}                  \label{2} \\
  1-s < \tfrac 2p                     \label{3} \\
  s+\tfrac{2}{p}+\tfrac{2}{r}   <3    \label{4} \\
  2 < s+\tfrac{2}{p}+\tfrac{2}{r}     \label{5} \\
  3 < s+\tfrac{2}{p}+\tfrac{4}{r},     \label{6}
\end{eqnarray}
then the map 
$\mathcal E \ni u \mapsto B(u) \in  \mathcal F^{\sscri 1}$ is well
defined and bounded. Therefore $\Phi: \mathcal E \to
\mathcal F$ is analytic. Moreover, there exist    positive
constants $C_1$ and $\ep$ such that
$$
  \|B\| \leq \;C_1\; T\,^\ep.
$$
\item[iii)]
Under the same assumptions as in $ii)$,
for any forcing term $f \in L^r(0,T;B^{-s}_{p\,q})$
and initial data $u_0 \in B^{-s+2-\frac{2}{r}}_{p\,r}$
there exists a time interval $[0,\overline T] \subseteq [0,T]$  and
a unique solution $u$ to equation \eqref{ns} with initial velocity
$u_0$, defined on the
time interval $[0,\overline T]$ with $\overline T>0$ satisfying
\eqref{picc}, and
\begin{align*}
& u \in  L^r(0,\overline T;B^{-s+2}_{p\,q}) \cap C([0,\overline T];
   B^{-s+2-\frac{2}{r}}_{p\,r})
\\
& u^\prime \in L^r(0,\overline T;B^{-s}_{p\,q})
\end{align*}
The solution $u$ depends analytically on the data $u_0$ and $f$.
\end{itemize}
\end{theorem}

The properties $i)$ of the linear operator  $d_0 \Phi$
follow from a general result:
\begin{proposition} \label{isom}
Let $T \in (0,\infty]$, $1<p,q,r<\infty$ and $s \in \rr$.\\
For any $f \in L^r(0,T;B^{-s}_{p\,q})$ and
$u_0 \in B^{-s+2-\frac{2}{r}}_{p\,r}$,
 there exists a unique $u \in W^{1,r}(0,T)\equiv
\{ u \in L^r(0,T;B^{-s+2}_{p\,q}): u^\prime \in L^r(0,T;B^{-s}_{p\,q})\}$
such that
\begin{equation*}
\begin{gathered}
 \Bigg\{
\end{gathered}
 \begin{aligned}
  u^\prime(t)+Au(t)&=f(t), \qquad t \in (0,T]\\
  u(0)&=u_0
 \end{aligned}
\end{equation*}
Moreover, the functions $u^\prime,u$ depend continuously on the data
$f$ and $u_0$,
that is there exists a positive constant $c$ such that
$$
\Big( \textstyle \int_0^T ( \|u(t)\|^r_{B^{-s+2}_{p\,q}}
+  \|u^\prime(t)\|^r_{B^{-s}_{p\,q}})\,dt\Big)^{1/r}
\leq
\Big( c \textstyle\int_0^T \|f(t)\|^r_{B^{-s}_{p\,q}} dt\Big)^{1/r} +
\|u_0\|_{B^{-s+2-\frac{2}{r}}_{p\,r}}
$$
Finally, the space $W^{1,r}(0,T)$
is continuously embedded in the space
$C_b([0,T];B^{-s+2-\frac{2}{r}}_{p\,r})$,
that is there exists
a positive constant $c$ such that
\begin{equation}\label{int-cont}
 \|u\|_{C_b([0,T];B^{-s+2-\frac{2}{r}}_{p\,r})}  \leq \;c\; \|u\|_{W^{1,r}(0,T)}
\end{equation}
and therefore the initial condition makes sense.\\
All the constants \footnote{We make the convention to denote
  different constants by the same symbol $c$, unless we want to mark them
  for further reference.}  depend only on $r,p,q,s$.
\end{proposition}
For the proof, see, e.g., \cite{b} Proposition 4.1
(based on \cite{dv}).
The assumptions on the linear operator $A$ and on the space $B^{-s}_{p\,q}$
are fulfilled; namely, the properties of the Stokes operator recalled
in Section 2 and the fact that $B^{-s}_{p\,q}$
($1<p,q<\infty$, $s \in \rr$) is a UMD Banach space (i.e. has the
Unconditional Martingale Difference property).
The last part of Proposition \ref{isom} is obtained by interpolation,
bearing in mind the interpolation result for Besov spaces
$B^{-s+2-\frac{2}{r}}_{p\,r}
           =(B^{-s}_{p\,q},B^{-s+2}_{p\,q})_{1-\frac{1}{r},r}$
for $1<r<\infty$, $1\leq p,q \leq \infty$, $s \in \rr$.

This corresponds to the first part $i)$ stated in Theorem \ref{t1}.
Assuming $i)$ and $ii)$, then part $iii)$ is proven as described
at the beginning of this section.
Hence, the proof of Theorem \ref{t1}  is complete as soon as we
prove $ii)$. This is given by the following result for the nonlinear equation.
\begin{proposition}   \label{quad}
Suppose the real numbers $s$ and $p,q,r \in (1,\infty)$ satisfy
the  conditions  (\ref{1}-\ref{6}) of Theorem \ref{t1}. Then
there exist constants $\ep>0$ and $C_1>0$ such that for all $T>0$
$$
  \left( \int_0^T \|B(u(t))\|^r_{B^{-s}_{p\,q}} dt \right)^{1/r}
 \leq
  \; C_1 \; T\,^\ep\; \|u\|^2_{W^{1,r}(0,T)},
  \qquad u \in W^{1,r}(0,T).
$$
\end{proposition}
{\bf Proof.}
First, let us show that there exists a pair $(a,b)$ of real numbers such that
\begin{eqnarray}
2-\tfrac 2r -s < a,b \label{U}\\
a,b < \tfrac 2p \label{D}\\
a+b=\tfrac 2p + 1-s \label{T}\\
\tfrac r2 \,[ (a+b)+2(s-2+\tfrac 2r)] < 1 \label{Q}
\end{eqnarray}
Before discussing these inequalities, let us remark that   \eqref{2} and \eqref{T} imply
that
\begin{equation} \label{C}
 a+b>0
\end{equation}
We begin from \eqref{5}, written as
$$
 2-\tfrac 2r -s < \tfrac 2p
$$
which grants that there are solutions to \eqref{U}-\eqref{D}.
\\
Secondly, \eqref{Q} is equivalent with
\begin{equation}\label{Qp}
 a+b < 4-2s-\tfrac 2r
        \tag{\ref{Q}'}
\end{equation}
Since by \eqref{4} (written as $\frac 2p +1-s < 4-2s-\frac 2r$)
any solution to \eqref{T} satisfies \eqref{Qp}
and hence \eqref{Q}, we only need to show that the system
\eqref{U}-\eqref{D}-\eqref{T} has at least one solution. We look
for a solution such that $a=b$. (Anyway, it is not difficult to
see that $a=b$ is not the only possible solution.) For this it is
enough that
$$
 2 -\tfrac 2r -s < \tfrac 1p +\tfrac 12 -\tfrac s2 < \tfrac 2p
$$
The second of this inequalities reads
$$
 1-s < \tfrac 2p
$$
which is \eqref{3}.
\\
The first one reads
$$
 2 - \tfrac 2r < \tfrac 1p + \tfrac 12 + \tfrac s2
$$
which is equivalent to \eqref{6}.
Thus, system (\ref{U}-\ref{Q}) has at least a solution.
\\
Define
$$
 \alpha = \tfrac r2 \,[a+s-2+\tfrac 2r]
$$
$$
 \beta = \tfrac r2 \, [b+s-2+\tfrac 2r]
$$
By \eqref{U} we have  $\alpha, \beta>0$ and  by \eqref{Q} we have
$\alpha+\beta<1$. In particular $\alpha,\beta<1$.
\\
We are now ready to finish the proof of Proposition
\ref{quad}. First we estimate the bilinear operator by means of
Bony's paraproducts techniques, as given in \cite{ch}, Corollary
1.3.1. Because $a$ and $b$ satisfy \eqref{D}-\eqref{T}-\eqref{C},
we have
\begin{equation} \label{sti}
\begin{split}
\|B(u(t))\|_{B^{-s}_{p\,q}} &=
 \|\nabla \cdot [u(t)\otimes u(t)]\|_{B^{-s}_{p\,q}}\\
 &\leq \; \|u(t)\otimes u(t)\|_{B^{-s+1}_{p\,q}}\\
 &\leq \;c\; \|u(t)\|_{B^a_{p\,q}} \|u(t)\|_{B^b_{p\,q}}\\
\end{split}
\end{equation}
Secondly, we use well known results on Besov spaces as interpolation
spaces (see, e.g., \cite{bl}) to get
$$
 \|u\|_{B^a_{p\,q}}
\leq
 \;c\;\|u\|^{1-\alpha}_{B^{-s+2-\frac{2}{r}}_{p\,q}}
   \|u\|^\alpha_{B^{-s+2}_{p\,q}}
$$
$$
 \|u\|_{B^b_{p\,q}}
\leq
   \;c\;\|u\|^{1-\beta}_{B^{-s+2-\frac{2}{r}}_{p\,q}} \|u\|^\beta_{B^{-s+2}_{p\,q}}
$$
with the interpolation parameters $\alpha, \beta$ defined above.
Here and in the following, $c$ denotes different constants.
\\
We use these inequalities to continue the estimate of the quadratic
operator from the last line of \eqref{sti}:
\begin{equation}\label{sti2}
\begin{split}
\|B(u(t))\|_{B^{-s}_{p\,q}}
&\leq \;c\; \|u(t)\|^{2-\alpha-\beta}_{B^{-s+2-\frac{2}{r}}_{p\,q}}
           \|u(t)\|^{\alpha+\beta}_{B^{-s+2}_{p\,q}}
\\
&\leq \;c\; \|u(t)\|^{2-\alpha-\beta}_{B^{-s+2-\frac{2}{r}}_{p\,r}}
           \|u(t)\|^{\alpha+\beta}_{B^{-s+2}_{p\,q}},
    \qquad \text{for }  r \leq q
\\
\end{split}
\end{equation}
Since $\alpha+\beta<1$, then
\begin{equation*}
 \begin{split}
  \big( \textstyle \int_0^T \|B(u(t))\|^r_{B^{-s}_{p\,q}} dt  \big)^{1/r}&
  \leq \;c\;
   \|u\|^{2-\alpha-\beta}_{C([0,T];B^{-s+2-\frac{2}{r}}_{p\,r})}
   \big( \textstyle\int_0^T  \|u(t)\|^{(\alpha+\beta)r}_{B^{-s+2}_{p\,q}}
        dt\big)^{1/r}
   \\
  &\leq \;c\;
    \|u\|^{2-\alpha-\beta}_{C([0,T];B^{-s+2-\frac{2}{r}}_{p\,r})}\;
      T^{(1-\alpha-\beta)/r}
    \big( \textstyle\int_0^T \|u(t)\|^r_{B^{-s+2}_{p\,q}}
  dt\big)^{(\alpha+\beta)/r}
   \\
  &\leq \;c\;
    \|u\|^2_{W^{1,r}(0,T)} T^{(1-\alpha-\beta)/r} \\
\end{split}
\end{equation*}
where we have used H\"older's inequality for the time integral in
the second line and the embedding \eqref{int-cont} of
Proposition \ref{isom} in the third line.
\hfill $\bsquare$

\vspace{1mm}
We point out that conditions \eqref{2} and \eqref{4}
are  against each  other;
namely,
rewritten down for the
regularity value of the initial data, they are
$$
 -s+2-\tfrac{2}{r} > -\tfrac{2}{r}-\big( \tfrac{2}{p}-1 \big)\qquad
   \mbox{ and } \qquad
 -s+2-\tfrac{2}{r} > \tfrac{2}{p}-1
$$
It follows that
$$
  -s+2-\tfrac{2}{r} >
  \max\Big\{ \big(-\tfrac 2r -(\tfrac 2p-1)\big), \; \big(\tfrac 2p -1\big) \Big\}
$$
By the computations \eqref{max2} in Appendix A, we get that
$ -s+2-\tfrac{2}{r} > -1$ at least. This imposes a restriction on the
admissible initial velocity to solve equation \eqref{ns} locally in time.

\section{Global existence and uniqueness} \label{sec:globale}
We want to show that the local solution constructed in the previous section
exists on the whole time interval $[0,T]$. To prove this, we split our
problem into two subproblems,
 considering two auxiliary variables $x$ and $y$ such that $u=x+y$.
\\
Following \cite{gp}, let decompose the data as
$$
 u_0=x_0+y_0 \; \text{ and }\; f=g+h
$$
The problem for the variable $y$ will have small forcing
term $h$ and small initial data $y_0$. Time-global existence and uniqueness will be
proved by means of Vishik and Fursikov's technique.
\\
The problem for the variable $x$ will
have more regular data: initial data $x_0 \in \mathcal D$ and force $g \in
L^2(0,T;H^{-1}_2)$. Time-global existence and uniqueness will be proved by
means of an  a priori estimate of the energy.

By the very definition of Besov spaces \eqref{def:b},
the space of periodic diver\-gence-free smooth
functions $\mathcal D$ is dense in any $B^s_{p\,q}$.
Therefore,
we have the following Lemma for the splitting of the data.
\begin{lemma} \label{splitta}
Let $f \in L^r(0,T;B^{\sigma}_{p\,q})$  with $1 \leq r,p,q<\infty$ and $\sigma \in \rr$.
Then for any $\ep>0$ there exist functions
$g^\ep \in
C^\infty([0,T];\mathcal D)$ and $h^\ep \in L^r(0,T;B^{\sigma}_{p\,q})$
such that $\| h^\ep\|_{L^r(0,T;B^{\sigma}_{p\,q})} < \ep$ and $f=g^\ep+h^\ep$.
\end{lemma}
The case of constant (in time) functions gives the splitting for the
initial data.

We proceed now in this way.
The two subproblems read
\begin{equation}\label{eq:x}
\left\{
  \begin{aligned}
   x^\prime(t)+Ax(t)+B(x(t),x(t))\qquad\qquad\qquad&\\
      +B(x(t),y(t))+B(y(t),x(t))&=g(t), \qquad t \in (0,T]\\[1mm]
   x(0) &=x_0
 \end{aligned}
\right.
\end{equation}
and
\begin{equation} \label{eq:y}
\begin{gathered}
 \Bigg\{
\end{gathered}
 \begin{aligned}
 y^\prime(t)+Ay(t)+B(y(t),y(t))&=h(t), \qquad t \in (0,T]\\
 y(0)&=y_0
 \end{aligned}
\end{equation}
For the latter one, we can choose $h$ and $y_0$ sufficiently small
(written as $\lll 1$ below),
in order to have the following result.
\begin{proposition} \label{y-globale}
Let the assumptions (\ref{1}-\ref{6}) of Theorem \ref{t1} be
satisfied.\\
Then, given any $y_0 \in B^{-s+2-\frac{2}{r}}_{p\,r}$, $h \in L^r(0,T;B^{-s}_{p\,q})$  with
$ \|y_0\|_{B^{-s+2-\frac{2}{r}}_{p\,r}}\lll 1$, $\|
h\|_{L^r(0,T;B^{-s}_{p\,q})} \lll 1$,
there exists a unique  solution $y$ to equation \eqref{eq:y} on the
time interval $[0,T]$, such that
\begin{equation}\label{reg:y}
  \begin{split}
   &y \in L^r(0,T;B^{-s+2}_{p\,q}) \cap C([0,T];B^{-s+2-\frac{2}{r}}_{p\,r}),\\
   &y^\prime \in L^r(0,T;B^{-s}_{p\,q}).   \\
  \end{split}
\end{equation}
\end{proposition}
{\bf Proof.}  Global existence  for small initial data and small
forcing term is obtained as described at the beginning of section \ref{sec:locale}
by means of Vishik and Fursikov's technique. Therefore,
parts $i)$ and $ii)$ of Theorem \ref{t1} entails this Proposition.
\hfill \bsquare

\vspace{1mm}

Concerning equation \eqref{eq:x}, we would like to show existence and
uniqueness on the whole  time interval $[0,T]$. We already know that
there exist $\overline T \in (0,T]$ and  a  function $x \,(=u-y)
\in  L^r(0,\overline T;B^{-s+2}_{p\,q})
\cap C([0,\overline T]; B^{-s+2-\frac{2}{r}}_{p\,r})$, $
x^\prime \in L^r(0,\overline T;B^{-s}_{p\,q})$.
To show that $x$ is indeed a solution to \eqref{eq:x}
is enough to show that all the terms in \eqref{eq:x} make sense. This is done,
for the nonlinear terms, analogously as in \eqref{sti}-\eqref{sti2}. For instance
$$
 \|B(x,y)\|_{B^{-s}_{p\,q}}\leq
           \;c\; \|x\|_{B^a_{p\,q}} \|y\|_{B^b_{p\,q}}
       \leq
  \|x\|^{1-\alpha}_{B^{-s+2-\frac{2}{r}}_{p\,q}} \|x\|^{\alpha}_{B^{-s+2}_{p\,q}}
   \|y\|^{1-\beta}_{B^{-s+2-\frac{2}{r}}_{p\,q}} \|y\|^{\beta}_{B^{-s+2}_{p\,q}}
$$
Condition \eqref{4} says that $ \dfrac 2p + \dfrac 2r +s<3$ (i.e.
$\dfrac 2p+\dfrac 2r -1<-s+2$ and $\dfrac 2p -1<-s+2-\dfrac
2r$).  Hence
  $B^{-s+2}_{p\,q} \subset B^{\frac 2p+\frac 2r -1}_{p\,q}$ and
  $B^{-s+2-\frac{2}{r}}_{p\,r} \subset B^{\frac 2p -1}_{p\,r}$
  and therefore
\begin{equation}\label{stimex1}
 \begin{split}
 &\exists \overline T \in (0,T] \mbox{ and } x \mbox{ such that }\\
 &\qquad x \in  L^r(0,\overline T;B^{\frac 2p+\frac 2r -1}_{p\,q})
 \cap C([0,\overline T];B^{\frac 2p -1}_{p\,r}),\\
 &\qquad x^\prime \in L^r(0,\overline T;B^{\frac 2p+\frac 2r -3}_{p\,q})  \\
 \end{split}
\end{equation}
Now we look for a priori estimates.
By \cite{gp} (Lemma 1.1), we  have the following energy estimate
\begin{lemma} \label{energ}
 Let $2 \leq p <\infty, 2< q <\infty, \dfrac 2p+\dfrac 2q -1>0$.
 Then, for any $\ep >0$ there exists a constant $c_\ep>0$ such that
 $$
  \left|\int_0^T \langle B(x(t)),y(t) \rangle \, d t \,\right|
  \leq
   \ep \int_0^T \|x(t)\|^2_{H^1_2} dt +
   c_\ep \, \int_0^T \|x(t)\|^2_{H^0_2}
   \|y(t)\|^q_{B_{p\,q}^{\frac{2}{p}+\frac{2}{q}-1}} dt.
 $$
\end{lemma}
A precisation on the proof is required. In fact, Gallagher
and Planchon in \cite{gp} work in the whole space. Anyway, the
technique used by them for the spatial domain $\mathbb D$
is $\rr^2$, works also in our case $\mathbb D=\toro$. Indeed, let
$u =\sum_k u_k e_k$ with $e_k=e_k(\xi)$   defined for $\xi \in
\rr^2$; the Littlewood--Paley decomposition (see, e.g.,
\cite{ch-libro} and references therein) gives
$$
 \Delta_m u= \sum_{2^m<|k|\leq 2^{m+1}} u_k e_k 
$$
where
$\Delta_m$ means the convolution with a function $\psi_m$ whose Fourier transform
$\hat \psi_m$ has support in
$\{ \xi \in \rr^2: 2^{m} <|\xi| \leq  2^{m+1}\}$.
Hence, the proof of \cite{gp}  based on Bony's paraproduct and
Ber{\v{n}}ste{\u{\i}}n's inequality\footnote{See, e.g.,
  \cite{ch-libro} for the definition and  properties of Bony's
  paraproduct and
  \cite{n} for Ber{\v{n}}ste{\u{\i}}n's  inequality.},
is valid also if we deal with
the Besov spaces
$$
  B^s_{p\,q}(\toro)= \{ u \underset{\mathcal{D}^\prime}{=} \textstyle\sum_{k\not 0} u_k e_k \in \mathcal{D}^\prime:
                         \textstyle\sum_{m \in \mathbb N}  \big( 2^{ms}
                         \|\Delta_m u\|_{L_p(\toro)} \big)^q<\infty
                     \}
$$
and  the $L_p(\toro)$-norms appear at the place of the
$L_p(\rr^2)$-norms of \cite{gp}. Moreover,  the norms 
$\|\nabla u\|_{L_2(\toro)} $ and $  \| u \|_{H^1_2(\toro)}$
are equivalent. Indeed,  $\toro$ is a bounded domain
and $\langle u, 1\rangle=0$ for $u\in \mathcal{D}^\prime$.

                  \vspace{1mm}
We now look for a priori estimates for the unknown $x$. Let us
multiply both sides of the first equation \eqref{eq:x} by $x$ and
integrate in space and in time. Two terms vanish, namely $\langle
B(y(t),x(t)),x(t)\rangle=0$ and $\langle
B(x(t),x(t)),x(t)\rangle=0$, see e.g. \cite{temam}.
Moreover $\langle B(x(t),y(t)),x(t)\rangle=-\langle
B(x(t),x(t)),y(t)\rangle$ and $\langle  A x(t),x(t)\rangle = \|
x(t) \|^2_{H^1_2}$. We then have
\begin{multline} \label{energia}
\tfrac 12 \|x(T)\|^2_{H^0_2} +
  \int_0^T \!\!\|x(t)\|^2_{H^1_2} \,dt
  = \tfrac 12 \|x_0\|^2_{H^0_2} + \int_0^T \langle B(x(t)),y(t) \rangle \, d t
  +\int_0^T \langle g(t),x(t)\rangle\,dt
\\
\underset{ \textrm{(by Lemma \ref{energ})} }{\leq}
  \tfrac 12 \|x_0\|^2_{H^0_2}
  +c \,
  \int_0^T \|x(t)\|^2_{H^0_2} \|y(t)\|^{\tilde q}_{B_{\tilde
           p\,\tilde q}^{\frac 2{\tilde p}+\frac 2{\tilde q} -1}}dt
\\
+ \tfrac{1}{2} \int_0^T \|x(t)\|^2_{H^1_2} \,dt
  +c \,\int_0^T\|g(t)\|^2_{H^{-1}_2} \,dt
\end{multline}
In this way, the required bounds are obtained by means of
Gronwall's lemma, as soon as we can find  $\tilde p \geq 2,
\tilde q >2$ with $\frac 2{\tilde p}+\frac 2{\tilde q}-1>0$, such
that 

$$
  \int_0^T \|y(t)\|^{\tilde q}_{B_{\tilde
            p\,\tilde q}^{\frac 2{\tilde p}+\frac 2{\tilde q}
            -1}}dt \; < \; \infty,
$$
where $y$ is the solution to problem \eqref{eq:y}. Proposition
\ref{y-globale} provides $y \in L^r(0,T;B^{-s+2}_{p\,q})$. Thus we
need to show that
\begin{equation} \label{regol:y}
   L^r(0,T;B^{-s+2}_{p\,q})
  \subseteq
   L^{\tilde q}(0,T;B_{\tilde p\,\tilde q}^{\frac 2{\tilde p}+\frac{2}{\tilde q}-1})
\end{equation}
for some $ \tilde p \geq 2, \tilde q >2$ with $\frac 2{\tilde p}+\frac
2{\tilde q}-1>0$.

Further conditions on the parameters $r,p,s$
are required in order that \eqref{regol:y}
holds.\\
First, there is the embedding
$L^r(0,T)\subseteq L^{\tilde q}(0,T)$ if
\begin{equation} \label{neltempo}
r \geq \tilde q
\end{equation}
since the time interval is finite.
\\
On the other hand, the space embedding
$B^{-s+2}_{p\,q} \subseteq B_{\tilde p\,\tilde q}^{\frac 2{\tilde
    p}+\frac{2}{\tilde q}-1}$
holds if
\begin{equation} \label{cond:utile}
 -s+2-\tfrac{2}{p} \geq \tfrac{2}{\tilde q}-1
\end{equation}
$$
 1  \leq q \leq \tilde q \leq \infty,
             \qquad
 1  \leq p \leq \tilde p \leq \infty
$$
(see \cite{bl} Theorem 6.5.1).
\\
We recall  assumption \eqref{4} of Theorem \ref{t1}:
$$
 s +\tfrac{2}{p} + \tfrac{2}{r} < 3
$$
The choice $\tilde q=r$, in order to satisfy \eqref{neltempo}, makes that the
conditions \eqref{4} and \eqref{cond:utile} become
identical (to be precise, there is the difference $<$ or $\leq $, which makes
\eqref{4} a slightly stronger than \eqref{cond:utile}).
Thus we have
\begin{equation}\label{nuova}
   \|y\|_{B_{p\,r}^{\frac 2p+\frac 2r -1}}
  \leq \;c\;
   \|y\|_{B^{-s+2}_{p\,q}}
\end{equation}

Summing up, choosing $\tilde q=r$ and $\tilde p=p$ in \eqref{energia} and
assuming (\ref{1}-\ref{6}) with with the additional conditions
$r>2, \frac 2p+\frac 2r -1>0 $,
the proper estimates follow. More precisely,
from \eqref{energia} we first have that
$$
 \|x(T)\|^2_{H^0_2}
 \leq
  \|x_0\|^2_{H^0_2}
  +c\, \int_0^T \|x(t)\|^2_{H^0_2}
             \|y(t)\|^r_{B^{-s+2}_{p\,q}}dt
  +c\,\int_0^T\|g(t)\|^2_{H^{-1}_2} \,dt
$$
Gronwall's lemma gives
$$
  \sup_{0\leq t \leq T} \|x(t)\|^2_{H^0_2}
  \le C(\|x_0\|_{H^0_2}, \|y\|_{L^r(0,T;B^{-s+2}_{p\,q})}, \|g\|_{L^2(0,T;H^{-1}_2)}) < \infty
$$
The last  result in conjunction with \eqref{energia} gives
$$
  \int_0^T \|x(t)\|^2_{H^1_2} \,dt
  \le C(\|x_0\|_{H^0_2}, \|y\|_{L^r(0,T;B^{-s+2}_{p\,q})}, \|g\|_{L^2(0,T;H^{-1}_2)} ) < \infty
$$
Finally
\begin{equation}\label{priori}
  \|x\|_{L^\infty(0,T;H^0_2)}+\|x\|_{L^2(0,T;H^1_2}
   \leq
  C(\|x_0\|_{H^0_2}, \|y\|_{L^r(0,T;B^{-s+2}_{p\,q})}, \|g\|_{L^2(0,T;H^{-1}_2)} )
\end{equation}

Use now the embedding theorem in Besov spaces
$$
 H^\sigma_2 \equiv B^\sigma_{2\,2}
       \underset{(p \geq 2)}\subseteq B^{\sigma-1+\frac 2p}_{p\,2}
       \underset{(\rho \geq 2)}\subseteq B^{\sigma-1+\frac 2p}_{p\,\rho}
$$
so to get from \eqref{priori} that
$$
  \|x\|_{L^\infty(0,T;B^{\frac 2p -1}_{p\,r}) }
  +\|x\|_{L^2(0,T;B^{\frac 2p}_{p\,q})}
   \leq
  C(\|x_0\|_{H^0_2}, \|y\|_{L^r(0,T;B^{-s+2}_{p\,q})}, \|g\|_{L^2(0,T;H^{-1}_2)} )
$$
By (complex) interpolation
  (that is $[L^2,L^\infty]_{1-\frac 2r}=L^r$ and
    $[B^{\frac 2p}_{p\,q},B^{\frac 2p -1}_{p\,q}]_{1-\frac 2r}
                   =B^{\frac 2p +\frac 2r -1}_{p\,q}$
    for $2<r<\infty$),
we obtain  that
\begin{equation} \label{reg:x}
  \|x\|_{L^\infty(0,T;B^{\frac 2p -1}_{p\,r})}
  +\|x\|_{L^r(0,T;B^{\frac 2p + \frac 2r -1}_{p\,q})}
   \leq
  C(\|x_0\|_{H^0_2}, \|y\|_{L^r(0,T;B^{-s+2}_{p\,q})}, \|g\|_{L^2(0,T;H^{-1}_2)} )
\end{equation}   
Comparing \eqref{reg:x} with \eqref{stimex1}, we get that the solution
$x$ exists on the whole time interval $[0,T]$.
We have therefore proven the following result.
\begin{proposition} \label{x-globale}
Let the assumptions (\ref{1}-\ref{6}) of Theorem \ref{t1} be
satisfied and moreover assume $r>2$ and $\frac 2p+\frac 2r -1>0$.
Then given any $x_0 \in \mathcal D$ and $g \in L^2(0,T;H^{-1}_2)$,
there exists a unique solution $x$ to equation \eqref{eq:x} on the
time interval $[0,T]$ such that
\begin{equation*}
 \begin{split}
 &x \in L^r(0,T;B^{\frac 2p+\frac 2r -1}_{p\,q})
 \cap C([0,T];B^{\frac 2p -1}_{p\,r})\\
 &x^\prime \in L^r(0,T;B^{\frac 2p+\frac 2r -3}_{p\,q})
 \end{split}
\end{equation*}
\end{proposition}
{\tt Remark.} Notice that more regularity on $g$ does not improve the
regularity of $x$, because of the presence of $y$ (which has a role similar to
an external force in equation \eqref{eq:x}). For this reason, we
assume $g \in  L^2(0,T;H^{-1}_2)$ instead of the other possible choice
$g \in C^\infty([0,T];\mathcal D)$.
On the other hand, the initial data $x_0$ is chosen very smooth in
order to consider without problems  the continuity in time in the next
results.
\hfill $\Box$
\vspace{2mm}

We combine Proposition \ref{y-globale} and Proposition
\ref{x-globale}
and, bearing in mind the embedding used to show \eqref{stimex1} (that is
to show that $x$ is less regular than $y$), we get that
\begin{equation}\label{reg:u}
 \begin{split}
 &u = x    +y \in
     L^r(0,T;B^{\frac 2p+\frac 2r -1}_{p\,q}) \cap C([0,T];B^{\frac 2p -1}_{p\,r}) \\
 &u^\prime = x^\prime +y^\prime \in
                  L^r(0,T;B^{\frac 2p+\frac 2r -3}_{p\,q})      \\
\end{split}
\end{equation}
This implies that there exists a function $u$, given by $u=x+y$, with the regularity specified
in \eqref{reg:u}. This is indeed a solution to equation \eqref{ns}.
In fact, bearing in mind equations \eqref{eq:x} and \eqref{eq:y}, we notice that
the function $u=x+y$ solves the equation \eqref{ns}, thanks to the fact that
the nonlinearity $B(u(t))$ is well defined; and this is so, because Chemin's result
to estimate the quadratic term garantees that this exists if
$u$ belongs to some Besov space of positive index and  from \eqref{reg:u} we have
that $u(t) \in B^{\frac 2p+\frac 2r -1}_{p\,q}$ with $\frac 2p+\frac 2r -1>0$
(for a.e. $t$).
Hence, $u$ is the sought global solution.

We sum up all the results proven so far and state  our main theorem.
\begin{theorem} \label{teo-gen}
For any forcing term $f \in L^r(0,T;B^{-s}_{p\,q})$
and initial velocity $u_0 \in B^{-s+2-\frac{2}{r}}_{p\,r}$
with
\begin{equation*}
\begin{split}
  s \in \rr, \qquad 1<p,q<\infty     \\
  2< r \leq q                        \\
  \tfrac 2p+\tfrac 2r -1>0           \\
  -\tfrac 2p <s -1< \tfrac{2}{p}                 \\
  2<s+\tfrac{2}{p}+\tfrac{2}{r}   <3  \\
  3 < s+\tfrac{2}{p}+\tfrac{4}{r}   \\
\end{split}
\end{equation*}
there exists a unique solution $u$ to equation \eqref{ns} on the
time interval $[0,T]$ such that
\begin{equation*}
\begin{split}
& u \in L^r(0,T;B^{\frac 2p+\frac 2r -1}_{p\,q}) \cap C([0,T];B^{\frac 2p -1}_{p\,r}) \\
& u^\prime  \in  L^r(0,T;B^{\frac 2p+\frac 2r -3}_{p\,q})\\
\end{split}
\end{equation*}
Moreover,
there exists a (strictly) positive $\overline T \le T$ such that
the above solution belongs to
$
L^2(0,\overline T;B^{-s+2}_{p\,q}) \cap C([0,\overline T];B^{-s+2-\frac{2}{r}}_{p\,r})
$.
Hence strong continuity for $t \to 0$ holds.
\end{theorem}
{\tt Remark}.
In  Appendix B, examples fulfilling
all these assumptions will be given.
The assumption $r>2$ imposes that $-s+2-\frac 2r >-\frac 12$,
as shown by \eqref{max3} in Appendix A.
\hfill $\Box$
\\[2mm]
{\bf Proof}. What remains to be proven is the uniqueness result.
Let us denote by  $\mathcal U$ the set of functions
$u$ satisfying the conditions \eqref{reg:u}. Consider two
solutions $u, \tilde u \in \mathcal U$ and denote by $\delta$ the
difference $u-\tilde u$. Then $\delta \in \mathcal U$ and it
satisfies the equation
\begin{equation} \label{diff}
 \begin{cases}
   \delta^\prime(t)+A\delta(t)+B(u(t),\delta(t))+B(\delta(t),\tilde u(t))
               =0, & \qquad t \in (0,T]\\
   \delta(0)=0&
 \end{cases}
\end{equation}
This is a linear equation in $\delta$.
We analyze this equation as a linear Stokes problem with a (linear)
perturbation term.
If the perturbation
$B(u,\cdot)+B(\cdot,\tilde u)$  is good enough (mainly, small for small time $T$,
so that this gives a small perturbation of the well-posed linear parabolic equation),
then there exists a unique solution. This will hold on a small time interval;
but since $\delta\equiv 0$ is a solution, then we get that the unique solution
is the zero one on a small time interval. Starting again from the zero value,
we can proceed in the same way to cover  the whole time interval. \\
We want to analyze the perturbation $B(u,\delta)+B(\delta,\tilde u)$.
We define the operators
$$
 \Gamma_u \delta:= - B(u,\delta) \qquad \tilde\Gamma_u \delta:= - B(\delta,u)
$$
and the space
$$
 S:=\{ \, \delta \in L^r(0,T;B^{\frac 2p+\frac 2r -1}_{p\,q}):
            \delta^\prime \in L^r(0,T;B^{\frac 2p+\frac 2r -3}_{p\,q})
                       \, \}
$$
equipped with the norm
$\|\delta\|_{S}=
     \big( \int_0^T \|\delta(t)\|^r_{B^{\frac 2p+\frac 2r -1}_{p\,q}} dt\;
      +\int_0^T \|\delta^\prime(t)\|^r_{B^{\frac 2p+\frac 2r -3}_{p\,q}} dt\big)^{1/r}
$.
\\
It is enough to consider the case with $u$,
since the same works for $\tilde u$, because of the symmetry of Chemin's estimates
in the two arguments.
We are going to show that
$$
 \Gamma_u: S \to L^r(0,T;B^{\frac 2p+\frac 2r -3}_{p\,q})
$$
and
\begin{equation}\label{domina}
 \| \Gamma_u \delta\|_{L^r(0,T;B^{\frac 2p+\frac 2r -3}_{p\,q})}
  \leq
  C_u(T) \; \|\delta \|_{S}
\end{equation}
with $C_u(T) \to 0$ as $T \to 0$.\\
This in nothing but an application of Chemin's estimates.
In fact
\begin{equation*}
\begin{split}
 \|B(u,\delta)\|_{B^{\frac 2p+\frac 2r -3}_{p\,q}}
 &\leq \|u \otimes \delta\|_{B^{\frac 2p+\frac 2r -2}_{p\,q}}
   \\
 &\leq  c\,
  \|u\|_{B^{\frac 2p+\frac 2r -1}_{p\,q}} \|\delta\|_{B^{\frac 2p -1}_{p\,q}}
  \qquad \textrm{ for } \tfrac 2p + \tfrac 2r -1 >0
   \\
 &\leq  c\,
  \|u\|_{B^{\frac 2p+\frac 2r -1}_{p\,q}} \|\delta\|_{B^{\frac 2p -1}_{p\,r}}
   \qquad  \textrm{ for } r \leq q \\
\end{split}
\end{equation*}
Therefore,
integrating in time, we get
\begin{multline}
  \|B(u,\delta)\|_{L^r(0,T;B^{\frac 2p+\frac 2r -3}_{p\,q})}
\\
 \leq   c
   \|u\|_{L^r(0,T;B^{\frac 2p+\frac 2r -1}_{p\,q})}
       \|\delta\|_{C([0,T];B^{\frac 2p -1}_{p\,r})}
\\
  \leq
    C_2 \|u\|_{L^r(0,T;B^{\frac 2p+\frac 2r -1}_{p\,q})}
    \|\delta\|_{S}.
\end{multline}
In the last step, we have used the interpolation result
(as \eqref{int-cont}) which allows to dominate the norm in
$C([0,T]; B^{\frac 2p -1}_{p\,r})$ by the norm in $S$
\begin{equation}\label{controlla}
  \|u\|_{C_b([0,T];B^{\frac 2p -1}_{p\,r})}\; \leq \;C\;
  \|u\|_{S}.
\end{equation}
Finally,
\eqref{domina} holds with
$$
 C_u(T)=  C_2 \|u\|_{L^r(0,T;B^{\frac 2p+\frac 2r -1}_{p\,q})}.
$$
We remark  that $C_u(T) \to 0$ as $T \to 0$.

We shall show that the problem
\begin{equation}
 \begin{cases}
   \delta^\prime(t)+A\delta(t)=-\Gamma_u \delta(t)-\tilde\Gamma_{\tilde u} \delta(t),
           & \qquad t \in (0,T],\\
   \delta(0)=0&
 \end{cases}
\end{equation}
has a unique solution $\delta \in S$ on a small time interval, using a contraction theorem.
Let us denote by $\Upsilon$ the mapping giving the solution to the Stokes problem
$$
 \Upsilon: g \mapsto v
 \qquad \begin{cases}
              v^\prime(t)+Av(t)=g\\
              v(0)=0
        \end{cases}
$$
We know from Proposition \ref{isom} that $\Upsilon$ is an
isomorphism from \\
 $L^r(0,T;B^{\frac 2p+\frac 2r
-3}_{p\,q})$ onto $S$ and
\begin{equation}\label{stokes}
 \| \Upsilon g \|_{S} \; \leq \; C_3\;
 \| g \|_{L^r(0,T;B^{\frac 2p+\frac 2r -3}_{p\,q})}
\end{equation}
for some positive constant $C_3$.
Therefore, by \eqref{domina} and \eqref{stokes} we get that
$$
 \|\Upsilon \circ (\Gamma_u \delta+\tilde \Gamma_{\tilde u} \delta) \|_{S}
\leq
 \; C_3 \,\big( C_u(T)+C_{\tilde u}(T) \big) \;
 \|\delta\|_{S}
$$
This shows that the mapping
$\Upsilon \circ (\Gamma_u+\tilde\Gamma_{\tilde u}): S \to S$ is a contraction
as soon as we work on the  time interval $[0,\overline T] \subseteq [0,T]$
with $\overline T>0$ chosen in such a way that
$$
 C_3 \,\big( C_u(\overline T)+C_{\tilde u}(\overline T) \big) \;<\; 1
$$
Hence there exists a unique solution $\delta \in S$ to equation
\eqref{diff} on the time interval $[0,\overline T]$.
This must coincide with the zero function: $\delta(t)=0$ for all $t \in [0,\overline T]$.
\\
Since the constants providing the contraction mapping depend only on the norms
of $u$ and $\tilde u$ (because problem \eqref{diff} is linear),
we start again from $\delta(\overline T)=0$ and we get the same result on the time
interval $[\overline T,2\overline T]$ and so on to conclude the proof
in a finite number of analogous steps.
\hfill $\bsquare$

\section{Remarks on the case where the forcing term $f$ is smooth}  \label{f-smooth}
In the previous section, the constraint $r>2$ has appeared.
This implies that the initial data $u_0$ is assumed to belong
 to the Besov space
$B^{-s+2-\frac{2}{r}}_{p\,r}$ with  $-s+2-\frac{2}{r}>-\frac 12$
(see  the comment at the end of section \ref{sec:locale}
and \eqref{max3} in Appendix A).
This restriction comes from the use of the energy estimate of Lemma \ref{energ}.
We want now to show that for $1<r\leq 2$,
when the  forcing term $f$ is smooth enough
( say, $f \in L^2(0,T;H^0_2)$ ),
a classical energy estimate can be used instead of Lemma \ref{energ}.
Therefore, assuming conditions (\ref{1}-\ref{6}) and some more regularity
on the forcing term, we prove global existence also when $1<r\leq 2$.
Notice that, when  the condition $r>2$ is removed, the initial velocity
regularity index $-s+2-\frac 2r$ can be very close to $-1$
(see also some examples of admissible values in  Appendix B).
This agrees with similar results obtained when there is no forcing term
(see, e.g., \cite{gp} for the problem in $\mathbb R^2$).
\\
We go back to the statement of Proposition \ref{y-globale} assuming
that $h=0$ (having taken $g=f$).
Since $y \in  L^r(0,T;B^{-s+2}_{p\,q})$, then
 there exists a $t_1 \in (0,T]$ as close to $0$ as we want, such that
$$
  y(t_1)\in B^{-s+2}_{p\,q}
$$
When the condition $1<r \leq 2$ is added to (\ref{1}-\ref{6}), then
$$
  -s+2 > 0
$$
because at the end of section \ref{sec:locale}
we have shown that $-s+2-\frac 2r >-1$.
\\
Moreover, if the index regularity of the initial velocity is negative, i.e.
$-s+2-\frac 2r<0$, then condition \eqref{4} imposes $p>2$.
At the end of Appendix B, it will be shown that this implies
$B^{-s+2-\frac 2r}_{p\,r}\not\subseteq B^0_{2\,2}$.
Therefore we are not dealing with the classical problem
of initial velocity with finite energy.
Summing up, if $1<r \leq 2, -s+2-\frac 2r<0$ and (\ref{1}-\ref{6}) hold, then
$$
  y(t_1)\in B^{-s+2}_{p\,2} \cap B^{ -s+2-\frac 2r}_{p\,r}
$$
Since $ -s+2-\frac 2r<0<-s+2$, there exists $\theta \in (0,1)$
such that
$$
 B^0_{p\,2} = \big(B^{-s+2-\frac 2r}_{p\,r}, B^{-s+2}_{p\,2}\big)_{\theta,2}
$$
Therefore
$$
 y(t_1) \in B^0_{p\,2} \subset B^0_{2\,2} \equiv H^0_2
$$
because the spatial domain is bounded. (In fact, in the same way we can show that
$y$ is infinitely smooth in space and time on the time interval $(0,T]$,
because there is no forcing term for $y$. But we do not need this result.)\\
Therefore on the time interval $[t_1,T]$,
the classical Hilbert-space theory can be applied to get existence and uniqueness results
(see, e.g., \cite{temam}):
$$
 y \in L^2(t_1,T;H^1_2) \cap C([t_1,T];H^0_2)
$$
By interpolation between $L^2(t_1,T;H^1_2)$ and $ L^\infty(t_1,T;H^0_2)$,
we get
$$
 y \in L^4(t_1,T;H^{1/2}_2)
$$

We now assume $f \in L^2(0,T;B^0_{2\,q})$.
Since
\begin{equation*}
\begin{split}
 &B^0_{2\,q} \subset B^{-s}_{p\,q} \qquad \textrm{ if } p>2, -1>-s-\tfrac 2p
 \\
 &L^2(0,T) \subseteq L^r(0,T)  \quad \textrm{ if } 1<r \leq 2
 \\
\end{split}
\end{equation*}
Then in our setting we have that
$$
 f \in L^2(0,T;B^0_{2\,q}) \subset L^r(0,T;B^{-s}_{p\,q})
$$
Therefore local existence results are obtained by means of
Proposition \ref{quad}. We only need  an a priori estimate on the
time interval $[t_1,T]$ for the unknown $x$. This is easily
obtained from the following classical estimate on the trilinear
term
\begin{alignat*}{2}
\langle B(x),y\rangle  & \leq \;\;\;\|x\|_{L_4} \|\nabla x\|_{L_2} \|y\|_{L_4}
                                &&\qquad\textrm{by  H\"older inequality}\\
  & \leq \;c\;\|x\|_{H^{1/2}_2} \|x\|_{H^1_2} \|y\|_{H^{1/2}_2}
                                && \qquad\textrm{by Sobolev embedding}\\
         & \leq \;c\; \|x\|^{1/2}_{H^0_2} \|x\|^{3/2}_{H^1_2} \|y\|_{H^{1/2}_2}
                               &&\qquad \textrm{by interpolation}\\
         & \leq \;\ep \; \|x\|^2_{H^1_2} + c_\ep \|x\|^2_{H^0_2}
                                \|y\|^4_{H^{1/2}_2}
                               &&\qquad  \textrm{by Young inequality}
\end{alignat*}
Hence in the time interval $[t_1,T]$,
the unknown $x$ does not explode in the required norms
and therefore there exists a unique $x \in L^2(t_1,T;H^1_2)\cap
C([t_1,T];H^0_2)$, $x^\prime \in L^2(t_1,T;H^{-1}_2)$.
We remind that $t_1$ can be chosen close to $0$ as much as we want.
Since existence on any small time interval $[0,t_1]$
was already proven in section \ref{sec:locale},
this result implies the global existence.
Finally, $u \in L^2(t_1,T;H^1_2)\cap C([t_1,T];H^0_2)$ for any $0<t_1<T$.
And this solution $u$ is unique.
\\[2mm]
{\tt Remark}.
 If $u_0\in H^0_2$ and $f \in L^2(0,T;H^{-1}_2)$, then a classical result
 grants that there exists a unique solution $u \in C([0,T];H^0_2) \cap
 L^2(0,T;H^1_2)$, $u^\prime \in  L^2(0,T;H^{-1}_2)$
 (see, e.g., \cite{temam}).
 \hfill $\Box$

\appendix
\newcommand{\appsection}[1]{\let\oldthesection\thesection
  \renewcommand{\thesection}{Appendix \oldthesection}
  \section{#1}\let\thesection\oldthesection}
\appsection{Lower estimates on $-s+2-\frac 2r$}
Let $1<p,q,r<\infty$ be given. Then
$$
\max \left\{-\tfrac 2r -\tfrac 2p +1,\, \tfrac 2p-1  \right\}
=
\begin{cases}
            \frac 2p -1    &    \quad 1<p\leq 2,\; \forall r\\[1mm]
            \frac 2p -1    &    \quad p>2,\,      1<r\leq\frac p{p-2}\\[1mm]
     -\frac 2r - \frac 2p +1  & \quad p>2,\,      r>\frac p{p-2}
           \end{cases}
$$
Hence, it easily follows that
\beq  \label{max2}
   \inf_{\substack{1<p<\infty \\1<r<\infty}}
      \max \left\{-\tfrac 2r -\tfrac 2p +1,\; \tfrac 2p-1  \right\}
   = -1
\eeq
With some more (but elementary) work, we obtain
\beq \label{max3}
 \inf_{\substack{1<p<\infty \\2<r<\infty}}
      \max \left\{-\tfrac 2r -\tfrac 2p +1,\; \tfrac 2p-1  \right\}
   = -\frac 12
\eeq

\appsection{Admissible values for the parameters $s,p,r$ and
numerical examples} We want to show that system (\ref{2}-\ref{6})
has a non void set of solutions. These are the conditions
appearing in Theorem \ref{t1}, providing local existence.
Two more conditions  are required for the global existence of
Theorem \ref{teo-gen}, unless the forcing term is smooth enough
(see Section \ref{f-smooth}). We start analyzing the
less restrictive conditions (\ref{2}-\ref{6}).\\
Set $x=\frac 2p$ and $y=\frac 2r$.
Then the system of conditions is
\begin{equation} \label{valori}
\begin{cases}
 2-s<x+y<3-s\\
 3-s<x+2y\\
 s-1<x\\
 1-s < x\\
 0<x,y<2
\end{cases}
\end{equation}
Because of the range of values specified in the last line,
necessary conditions for the existence of a solution
to \eqref{valori} are
$$
\begin{cases}
2-s<4 \mbox{ and } 0<3-s\\
3-s<6\\
s-1<2\\
1-s<2
\end{cases}
$$
Hence, the admissible values for the parameter $s$ are
$-1<s<3$.
Moreover $x<3-s$ (from the first and the last line in \eqref{valori})
and $x>s-1$ (that is the third line in \eqref{valori})
imposes the further restriction: $s<2$.
\\
Summing up, the admissible values for the parameter $s$ are
$$
  -1<s<2
$$
\\
We distinguish two cases.
\\
$\bullet  -1<s<1$
\\
Since $s-1<1-s$, the third line in  \eqref{valori} can be neglected.
Representing the remaining conditions \eqref{valori} on the $(x,y)$-plane, it is
easy to see that there exist solutions. If we are interested in the
solutions satisfying also the condition $-1<-s+2-y<0$ (for the
regularity of the initial velocity),
then $s$ must be positive. We give examples of parameters satisfying
the above conditions.
\\
Examples:
\begin{alignat*}{6}
 \text{ for Pro. \ref{quad} }
 &\qquad s=\tfrac 9{10}
     &\quad r= \tfrac{20}{19}
          &\quad  p=12
                &&\qquad -s+2-\tfrac 2r = -\tfrac 45\\
 \text{ for Th. \ref{teo-gen} }
 &\qquad s=-\tfrac 9{10}
     &\quad r=\tfrac{100}{49}
          &\quad  p=\tfrac{40}{39}
                &&\qquad -s+2-\tfrac 2r = \tfrac {48}{25}
\end{alignat*}
We do not choose $q$ since this is the less significant
parameter to characterize a Besov space.
 Notice that $1<r\leq 2$ in the first case, providing global existence
for ``regular'' forcing term. The second case concerns
 positive index regularity $-s+2-\frac 2r$ for the initial velocity.
\\[1mm]
$\bullet \,1\leq s <2 $
\\
Since $s-1\geq 1-s$, the forth line in  \eqref{valori} can be neglected.
Again the graphic representation shows that there are solutions.
The condition $-1<-s+2-y<0$ (for the
regularity of the initial velocity)
requires $x<1$ (i.e. $p>2$). Notice that in this case
both conditions $y<1$ (i.e. $r>2$) and $-s+2-y<0$ (i.e. $-s+2-\frac 2r<0$)
can be fulfilled, providing parameters satisfying the assumptions of
Theorem \ref{teo-gen} with the initial velocity in a Besov space of
negative order (not allowed in the previous case).
\\
Examples:
\begin{alignat*}{6} 
\text{for Pro. \ref{quad}}
 &\qquad s=\tfrac{11}{10},
     &\quad r=\tfrac{8}7,
          &\quad p=\tfrac{40}{3},
                  &&\qquad -s+2-\tfrac 2r = -\tfrac{17}{20}\\
 \text{ for Th. \ref{teo-gen} }
 &\qquad s=\tfrac{149}{100},
     &\quad r=\tfrac{200}{99},
          &\quad  p=4,
                &&\qquad -s+2-\tfrac 2r = -\tfrac{48}{100}\\
 &\qquad s=\tfrac{11}{10},
     &\quad r=\tfrac{40}{19},
          &\quad  p=3,
                &&\qquad -s+2-\tfrac 2r = - \tfrac{1}{20}\\
 &\qquad s=\tfrac 43,
     &\quad r=3,
          &\quad  p=\tfrac 52,
                &&\qquad -s+2-\tfrac 2r = 0\\
 &\qquad s=\tfrac {19}{10},
        &\quad r=21,
              &\quad p=2,
                   &&\qquad -s+2-\tfrac 2r = \tfrac{1}{210}
\end{alignat*}
This is the only case providing global solutions with initial velocity $u_0$
in Besov space of negative order $-s+2-\tfrac 2r<0$ and force
$f \in L^r(0,T;B^{-s}_{p\,q})$.
\\[1mm]
{\tt Remark.}
We remark that $-s+2-\frac 2r<0$ and $s+\frac 2p+\frac 2r<3$
imposes $\frac 2p<1$, that is $p>2$. This implies that there
the embedding $B^{-s+2-\frac 2r}_{p\,r}\subset B^0_{2\,2}$ never holds.
Hence we really deal with a generalization of the classical result
for $u_0\in B^0_{2\,2}$.
An analogous statement holds for the forcing term:
$L^r(0,T;B^{-s}_{p\,q}) \not\subset L^2(0,T;B^{-1}_{2\,2})$, because
$B^{-s}_{p\,q}\not\subset B^{-1}_{2\,2}$ for
$p>2$ and $s>1$.
 \hfill $\Box$

\vspace{2mm} \noindent {\bf Acknowledgments.} The second author
wishes to thank the Department of Mathematics, Hull University,
for the warm hospitality. Financial support by the Alexander von
Humboldt Stiftung is gratefully acknowledged. This work was
partially supported by British Council/DAAD grant. 


\begin{thebibliography}{99999}

       \bibitem[AC]{ac}
Albeverio, S.; Cruzeiro, A.B.
 {Global flows with invariant (Gibbs) measures for Euler and
    Navier-Stokes two dimensional fluids.}
 Comm. Math. Phys.  {\bf 1990}, {\tt 129} (3), 431-444.

      \bibitem[ARH-K]{afhk}
Albeverio, S.; Ribeiro de Faria, M.; H{\o}egh--Krohn, R.
 {Stationary measures for the periodic Euler flow in two dimensions.}
 J. Statist. Phys. {\bf 1979}, {\tt 20} (6), 585-595.

      \bibitem[AFa]{af-sw}
Albeverio, S.; Ferrario, B.
 {Invariant measures of L\'evy-Khinchine type for 2D fluids.}
 Proceedings  of the Swansea 2002 Workshop
 ``Probabilistic Methods in Fluids'',
 Wales, UK, 14-19 April 2002. 
 Eds.: I.M. Davies, N. Jacob, A. Truman, O. Hassan, K. Morgan, N.P. Weatherill.
 World Scientific (2003), 130-143

      \bibitem[AFb]{af}
Albeverio, S.; Ferrario, B.
 {Uniqueness of solutions of the stochastic
        Navier--Stokes equation with invariant measure given by the enstrophy.}
Annals of Probability  {\bf 2004} {\tt 23} (2), 1632-1649 

\bibitem[AH-K]{ahk}
 Albeverio, S.; H{\o}egh--Krohn, R.
 {Stochastic flows with stationary distribution for two-dimensional
    inviscid fluids.}
 Stochastic Process. Appl. {\bf 1989}, {\tt 31} (1), 1-31.

        \bibitem[BG]{bg}
Biagioni, H.; Gramchev, T.
 {On the 2D Navier--Stokes equation with singular initial data and
 forcing term.}
 Mat. Contemp. {\bf 1996}, {\tt 10}, 1-20.

      \bibitem[B91]{b91}
 Brze\'zniak, Z.
 {On analytic dependence of solutions of Navier--Stokes
 equations with respect to exterior force and initial velocity.}
 Univ. Iagel. Acta Math. {\bf 1991}, {\tt 28}, 111-124.

    \bibitem[B95]{b}
Brze\'zniak. Z.
 {Stochastic partial differential equations in M-type 2 Banach spaces},
 Potential Anal. {\bf 1995} {\tt 4} (1), 1-45.

    \bibitem[BL]{bl}
Bergh, L.; L\"ofstr\"om. J.
 {\tt Interpolation Spaces. An introduc\-tion};
 Springer: Berlin-New York, 1976.

    \bibitem[Ca]{ca}
Cannone, M.
 {A generalization of a theorem by Kato on Navier--Stokes equations},
 Rev. Mat. Iberoamericana {\bf 1997}, {\tt 13} (3), 515-541.

    \bibitem[CP]{CP}
Cannone, M.; Planchon, F.
 {On the non stationary Navier-Stokes equations with an external force},
Adv. Diff. Eq. {\bf 1999} {\tt 4} (5), 697-730.

    \bibitem[Ch96]{ch}
Chemin, J.-Y.
 {About Navier--Stokes system.}
 Publication du Laboratoire d'Analyse Num\'erique {\bf 1996}, {\tt R 96023}.

    \bibitem[Ch98]{ch-libro}
Chemin,  J.-Y.
 {\tt Perfect incompressible fluids};
 Oxford University Press: New York, 1998.

    \bibitem[DPD]{dpd}
Da Prato, G.; Debussche, A.
 {2D-Navier--Stokes equations driven by a space--time white noise.}
 J. Funct. Anal. {\bf 2002}, {\tt 196}  (1), 180-210.

    \bibitem[DV]{dv}
Dore, G.;  Venni, A.
 {On the closedness of the sum of two closed operators.}
 Math. Z. {\bf 1987}, {\tt 196} (2), 189-201.

    \bibitem[GP]{gp}
Gallagher, I.;  Planchon, F.
 {On infinite energy solution to the Navier--Stokes equations:
  global 2D existence and 3D weak-strong uniqueness.}
 Arch. Rat. Mech. Anal. {\bf 2002}, {\tt 161} (4), 307-337.

    \bibitem[L]{l}
Lunardi, A.
 {\tt Interpolation Theory};
 Appunti Scuola Normale Supe\-rio\-re: Pisa, 1999.

    \bibitem[N]{n}
Nikol'ski\u{\i},  S.M.
 {\tt Approximation of functions of several vari\-ables and imbedding theorems}; 
 Springer: New York-Heidelberg, 1975.

    \bibitem[Te79]{temam}
Temam, R.
 {\tt Navier-Stokes Equations: theory and nu\-meri\-cal analysis},
 Reprint of the 1984 edition;
 AMS Chelsea Publishing: Providence RI, 2001.

     \bibitem[Te83]{te-p} 
Temam, R.
 {\tt Navier-Stokes equations and nonlinear func\-tional analysis},
 Second edition; CBMS-NSF Regional Conference Series in Applied Mathematics, 66;
 SIAM: Philadelphia PA, 1995.

    \bibitem[VF]{vf}
Vishik, M.J.; Fursikov, A.V.
 {\tt Mathematical Problems of Stat\-is\-ti\-cal Hydromechanics},
 (translated from the Russian);
 Kluwer: Dordrecht-Boston-London, 1988.

\end{thebibliography}
\end{document}